\documentclass[times]{my-iapress}

\usepackage{moreverb}

\usepackage[colorlinks,bookmarksopen,bookmarksnumbered,citecolor=red,urlcolor=red]{hyperref} 

\def\volumeyear{(2019)}
\def\volumenumber{7}
\def\volumemonth{No 2}
\usepackage{amsmath,amssymb}
\usepackage{graphicx}

\usepackage{subcaption}

\newtheorem{definition}{Definition}

\newcommand{\Dl}{{}^{\scriptscriptstyle C}_{\scriptscriptstyle 0}\!D^{\alpha}_{\scriptscriptstyle t}}
\newcommand{\DrRL}{{}_{\scriptscriptstyle t}D^{\alpha}_{\scriptscriptstyle t_f}}

\begin{document}

\runningheads{Optimal Control and Sensitivity Analysis of a Fractional Order TB Model}{S. Rosa and D. F. M. Torres}

\title{Optimal Control and Sensitivity Analysis of a Fractional Order TB Model}

\author{Silv\'erio Rosa\affil{1}, Delfim F. M. Torres \affil{2}$^,$\corrauth}

\address{\affilnum{1}Department of Mathematics and Instituto de Telecomunica\c{c}\~{o}es (IT),
University of Beira Interior, 6201-001 Covilh\~{a}, Portugal.\\ 
\affilnum{2}Center for Research and Development in Mathematics and Applications (CIDMA),
Department of Mathematics, University of Aveiro, 3810-193 Aveiro, Portugal.}

\corraddr{Delfim F. M. Torres (Email: delfim@ua.pt). 
Department of Mathematics, University of Aveiro, 3810-193 Aveiro, Portugal.}


\begin{abstract}
A Caputo fractional-order mathematical model for the transmission dynamics 
of tuberculosis (TB) was recently proposed in 
[Math. Model. Nat. Phenom. 13 (2018), no.~1, Art.~9]. 
Here, a sensitivity analysis of that model is done, showing the importance of accuracy 
of parameter values. A fractional optimal control (FOC) problem is then formulated 
and solved, with the rate of treatment as the control variable. Finally, a 
cost-effectiveness analysis is performed to assess the cost and the effectiveness 
of the control measures during the intervention, showing in which conditions 
FOC is useful with respect to classical (integer-order) optimal control. 
\end{abstract}


\maketitle


\section{Introduction}

Tuberculosis (TB) is an infectious disease caused by the bacterium 
\emph{Mycobacterium tuberculosis}, which is usually spread through
the air when people who have active TB in their lungs cough, spit, speak, or sneeze.
TB is one of the top ten causes of death worldwide, which justifies the interest of
researchers on the area: see, e.g., \cite{MR3562914,MR3101449,WHO:TB:report}. 
A good  survey on optimal models of TB is presented in \cite{MR3644001}. 
Delays are introduced in a TB model in \cite{MR3562914}, representing the time delay 
on the diagnosis and initiation of treatment of individuals with active TB infection. 
Optimal control strategies to minimize the cost of interventions, considering reinfection 
and post-exposure interventions, are investigated in \cite{MR3101449}. In \cite{rodrigues2014cost}, 
the potential of two post-exposure interventions, treatment of early latent TB individuals 
and prophylactic treatment/vaccination of persistent latent TB individuals, is investigated. 
A mathematical model for TB is studied in \cite{MR3388961} from the optimal control point of view, 
using a multiobjective optimization approach. For numerical simulations using TB real data 
from Angola, see \cite{MR2970904}. In spite of the numerous works on tuberculosis, 
the literature on fractional-order mathematical models for TB is still scarce. 
In \cite{SWEILAM2016271}, it is proposed a multi-strain TB model with variable-order
fractional derivatives and a numerical scheme to approximate the endemic solution, numerically,
is developed. More recently, in 2018, a Caputo type fractional-order mathematical model 
for the transmission dynamics of tuberculosis was proposed, based on the nonlinear differential 
system studied in \cite{YANG201079}, and its stability is investigated \cite{TB:frac:2018}.
Such model is here subject to fractional optimal control theory \cite{MR3673702,MR3673710}, 
and its usefulness to tackle a TB epidemic scenario discussed.

The paper is organized as follows. In Section \ref{sec:model}, we introduce 
the fractional-order TB model. The main results are then given in Section~\ref{sec:mresults}: 
sensitivity analysis of the fractional TB model (Section~\ref{sec:sensitivity});
fractional optimal control, cost-effectiveness and numerical simulations for the TB model 
(Sections~\ref{sec:optctrl} and \ref{sec:nresults}). We end with 
Section~\ref{sec:conclusions} of conclusions.

\section{Fractional-order Tuberculosis Model}
\label{sec:model}

In this section, we consider a Caputo fractional-order tuberculosis (TB) model 
due to \cite{TB:frac:2018}. The model describes the dynamics of a population 
that is susceptible to infection by the \emph{Mycobacterium tuberculosis} with
incomplete treatment. The population consists of four compartments: susceptible 
individuals ($S$); latent individuals ($L$), which have been infected but are 
not infectious and do not show symptoms of the disease; infectious individuals
($I$), which have active TB, may transmit the infection, but are not under 
treatment; and  infected individuals in treatment ($T$).

\begin{figure}[!htb]
\centering
\includegraphics[scale=0.86]{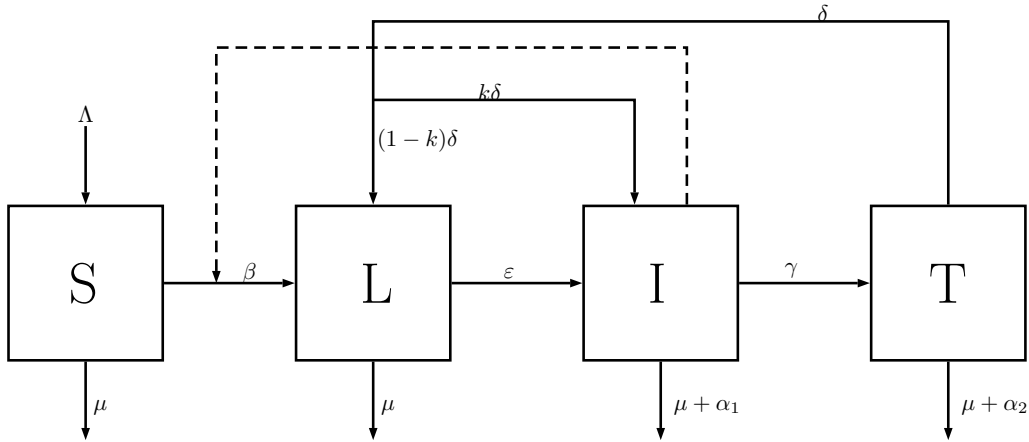}
\caption{Epidemiological scheme of the mathematical model \eqref{TB_model}.}
\label{f:diag}
\end{figure}

The susceptible population is increased by the recruitment of individuals 
into the population at a rate $\Lambda$. All individuals are exposed 
to natural death, at a constant rate $\mu$. Deterministic continuous transitions 
between the compartments, also known as states, are used. Susceptible individuals 
acquire TB infection by the contact with infected individuals at a rate $\beta I$, 
where parameter $\beta$ is the transmission coefficient. Individuals in the latent 
class, $L$, become infectious at a rate $\varepsilon$, and infectious individuals, 
$I$, start treatment at a rate $\gamma$. Treated individuals, $T$, leave their 
compartment at rate $\delta$. After leaving the treatment compartment, an individual 
may enter compartment $L$, due to the remainder of \emph{Mycobacterium tuberculosis}, 
or compartment $I$, due to the failure of treatment. The parameter $k$, 
$0 \leqslant k \leqslant 1$, represents the failure of the treatment, 
where $k = 0$ means that all the treated individuals shall become latent, 
while $k = 1$ means that the treatment fails and all the treated individuals 
shall still be infectious. Infectious, $I$, and under treatment individuals, $T$, 
may suffer TB-induced death at the rates $\alpha_1$ and $\alpha_2$, respectively. 
The Caputo fractional-order system of differential equations that describes the
above assumptions is
\begin{equation}
\label{TB_model}
\begin{cases}
\Dl S(t) =  \Lambda-\beta I(t)S(t)-\mu S(t),\\
\Dl L(t) = \beta I(t) S(t)+(1-k)\delta T(t)-(\mu+\varepsilon)L(t),\\
\Dl I(t) =\varepsilon L(t)+k\delta T(t)-(\mu+\gamma+\alpha_1)I(t),\\
\Dl T(t) =\gamma I(t) -(\mu+\delta+\alpha_2)T(t),
\end{cases}
\end{equation}
where  $\Dl$ denotes the left Caputo derivative of order $\alpha \in (0,1]$
\cite{MR1658022,MR3736617}. Note that when $\alpha=1$, the fractional 
compartmental model \eqref{TB_model} represents the classical TB model 
studied in \cite{YANG201079}. In Figure~\ref{f:diag}, one can find
a diagram representing the model dynamics.

\section{Main results}
\label{sec:mresults}

We begin by investigating the sensitiveness of the TB model \eqref{TB_model},
discussed in Section~\ref{sec:model}, with respect to
the variation of each one of its parameters.

\subsection{Sensitivity analysis}
\label{sec:sensitivity}

One of the most significant thresholds when studying infectious disease models 
is the basic reproduction number \cite{van2017reproduction}. The basic 
reproduction number for the TB model \eqref{TB_model} is given by 
\begin{equation}
\label{r0:TB_model}
R_0=\frac{\beta\varepsilon b_3 \Lambda}{\mu b_1 b_2 b_3 
-\mu\delta\gamma((1-k)\varepsilon+k b_1)}
\end{equation}
with $b_1=\mu+\varepsilon$, $b_2=\mu+\gamma+\alpha_1$ 
and $b_3=\mu+\delta+\alpha_2$ \cite{YANG201079}.

Now we perform a sensitivity analysis
for the endemic threshold \eqref{r0:TB_model}.
Such analysis tells us how important each parameter is to disease transmission.
This information is crucial not only for experimental design,
but also to data assimilation and reduction of complex models
\cite{powell2005sensitivity}. Sensitivity analysis is commonly used
to determine the robustness of model predictions to parameter values,
since there are usually errors in collected data and presumed parameter values.
It is used to discover parameters that have a high impact on the threshold $R_0$
and should be targeted by intervention strategies. More accurately,
sensitivity indices's allows us to measure the relative change
in a variable when parameter changes. For that purpose we use the normalized
forward sensitivity index of a variable, with respect to a given parameter,
which is defined as the ratio of the relative change in the variable
to the relative change in the parameter. If such variable is differentiable
with respect to the parameter, then the sensitivity index
is defined using partial derivatives.

\begin{definition}[See \cite{chitnis2008determining,rodrigues2013sensitivity}]
\label{def:sentInd}
The normalized forward sensitivity index of $R_0$, which is differentiable
with respect to a given parameter $p$, is defined by
$$
\Upsilon_p^{R_0}=\frac{\partial R_0}{\partial p}\frac{p}{R_0}.
$$
\end{definition}

The values of the sensitivity indices for the parameters values of Table~\ref{tab:param},
are presented in Table~\ref{tab:sensitivity}. 

\begin{table}[ht!]
\centering
\caption{Values of the models' parameters taken from \cite{TB:frac:2018,YANG201079}.}
\begin{tabular}{lll}
\toprule
Name & Description & Value\\
\midrule
$\Lambda$ & Recruitment rate & 792.8571\\
$\beta$ & Transmission coefficient & $0.0005$\\
$\mu$ & Natural death rate & 1/70\\
$k$ & Treatment failure rate & 0.15\\
$\delta$ & Rate at which treated individuals leave the $T$ compartment & 1.5\\
$\varepsilon$ & Rate at which latent individuals $L$ become infectious & 0.00368\\
$\alpha_1$ & TB-induced death rate for infectious individuals $I$ & 0.3\\
$\alpha_2$ & TB-induced death rate for under treatment  individuals $T$ & 0.05\\
\bottomrule
\end{tabular}
\label{tab:param}
\end{table}
\begin{table}[!htb]
\centering
\caption{Sensitivity of $R_0$ evaluated for the parameter values given
in Table~\ref{tab:param}.}\label{tab:sensitivity}
\begin{tabular}{cc@{\hspace*{5cm}}c}\toprule
Parameter && Sensitivity index \\[1mm] \midrule
$\mu$ && $-1.93223$\\[1mm]
$\varepsilon$ && $+0.911803$\\[1mm]
$\gamma$ && $-0.605532$\\[1mm]
$\alpha_1$ && $-0.376538$\\
$\delta$ && $0.0112215$\\
$\alpha_2$ && $-0.0872783$\\
$\Lambda$ && $+1$\\
$k$ && $+0.100487$\\
$\beta$ && $+1$\\ \bottomrule
\end{tabular}
\end{table}

\begin{figure}[!htb]
\centering
\begin{subfigure}[b]{0.46\textwidth}\centering
\includegraphics[scale=0.46]{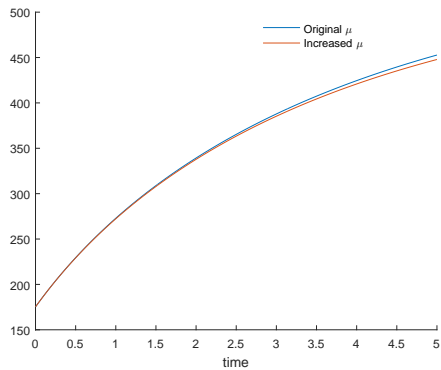}
\caption{Impact of the variation of $\mu$ in the number 
of infective individuals $I$
(there is a visible difference).}
\label{fig:sensit_mu}
\end{subfigure}\hspace*{1cm}
\begin{subfigure}[b]{0.46\textwidth}
\centering
\includegraphics[scale=0.46]{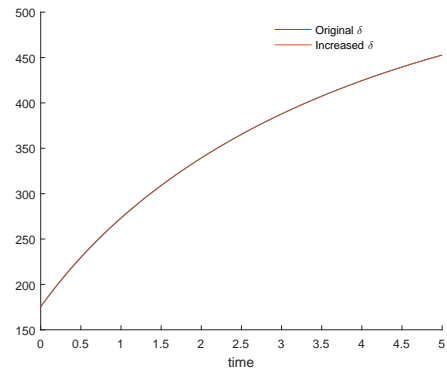}
\caption{Impact of the variation of $\delta$ in 
the number of infective individuals $I$ (no visible difference).}
\label{fig:sensit_delta}
\end{subfigure}
\caption{Infected number of individuals predicted by the TB model 
\eqref{TB_model} with original parameter values as in Table~\ref{tab:param}
and with an increase of 15\% of a specific parameter:
$\mu$ in (a) and $\delta$ in (b).}
\label{fig:sensitivity}
\end{figure}

Note that the sensitivity index may depend on several parameters 
of the system, but also can be constant, independent of any parameter. 
For example, $\Upsilon_{\Lambda}^{R_0}=+1$, meaning that increasing 
(decreasing) $\Lambda$ by a given percentage increases
(decreases) always $R_0$ by that same percentage.
The estimation of a sensitive parameter should be carefully done,
since a small perturbation in such parameter leads to relevant
quantitative changes. On the other hand, the estimation of
a parameter with a rather small value for the sensitivity index
does not require as much attention to estimate,
because a small perturbation in that parameter leads
to small changes \cite{mikucki2012sensitivity}.
According to our results in Table~\ref{tab:sensitivity}, 
special attention should be paid to the estimation of the 
death rate, $\mu$. In contrast, the rate at which individuals 
leave state $T$, $\delta$, does not require as much attention 
because of its low value of the sensitivity index.
This is well illustrated in Figure~\ref{fig:sensitivity},
where we can see in Figure~\ref{fig:sensit_mu}
the graphics of the number of infective with and without
an increment of 15\% for the parameter $\mu$ (most sensitive parameter),
and in Figure~\ref{fig:sensit_delta} with identical comparison 
of graphics for parameter $\delta$ (less sensitive parameter).

Using the values of the parameters as proposed in \cite{TB:frac:2018}, 
we obtain results very different from those that the authors present there. 
We noticed that the value of parameter $\mu$ presented in \cite{TB:frac:2018} 
($\mu'=0.143$), the most sensitive one, is erroneous, being 10 times the correct 
value \cite{YANG201079} and was not the value that authors' used in the simulation. 
This particular case highlights the importance of the analysis
we have carried out here.

\subsection{Fractional optimal control} 
\label{sec:optctrl}

The circumstances on which variation of the variables of the model depend, may be controlled.
In case of TB disease, the fraction of infectious individuals that is identified and put under
treatment is one of the most commonly used. Therefore, the parameter $\gamma$ is, in what 
follows, replaced by a control variable. Consequently, we consider the following fractional
optimal control problem: to minimize the number of infectious individuals and the cost
of the control of the disease with the identification and treatment of the patients, that is,
\begin{equation}
\label{cost-functional}
\min ~\mathcal{J}(I(t),u(t))
=\int_0^{t_f} I(t)+ B \rho\, u(t)^2 ~dt
\end{equation}
with $0<B <\infty$, subject to the fractional control system
\begin{equation}
\label{eq:modTB_control}
\begin{cases}
\Dl S(t) &=  \Lambda-\beta I(t)S(t)-\mu S(t)\\
\Dl L(t) &= \beta I(t) S(t)+(1-k)\delta T(t)-(\mu+\varepsilon)L(t)\\
\Dl I(t) &=\varepsilon L(t)+k\delta T(t)-(\mu+u(t)+\alpha_1)I(t)\\
\Dl T(t) &=u(t) I(t)-(\mu+\delta+\alpha_2)T(t)
\end{cases}
\end{equation}
and given initial conditions
\begin{equation}
\label{ocp:ic}
S(0),L(0),I(0),T(0)\geqslant 0.
\end{equation}
Here, $\rho$ is the  maximum number of infectious individuals 
of the problem without control; and the control variable, $u$, 
designate the fraction of individuals that is put under treatment. 
The set of admissible control functions is
\begin{equation}
\label{Omega:set}
\Omega=\left\{u(\cdot)\in L^{\infty}(0,t_f):
0\leqslant u(t)\leqslant u_{\max},\forall t\in[0,t_f]\right\}.
\end{equation}
Pontryagin's maximum principle (PMP) for fractional optimal control can be used
to solve the problem \cite{MR3529374,MR3443073,MR3225198,malinowska2012introduction}.
The Hamiltonian associated with our optimal control problem is
\begin{multline*}
\mathcal{H}
= I+ B \rho\, u^2 + p_1(\Lambda-\beta I S-\mu S)\\
+p_2(\beta I S+(1-k)\delta T-(\mu+\varepsilon)L)
+p_3(\varepsilon L+k\delta T-(\mu+u+\alpha_1)I)
+p_4(u I-(\mu+\delta+\alpha_2)T);
\end{multline*}
the adjoint system uphold that the co-state variables $p_i(t)$, 
$i=1,\ldots,4$, verify
\begin{equation}
\label{eq:co_states_fr2}
\begin{cases}
\DrRL p_1(t) = \mu p_1(t) - \beta I(t) (p_1(t)-p_2(t)),\\[1.2mm]
\DrRL p_2(t) =  (\varepsilon+\mu) p_2(t) - \varepsilon p_3(t), \\[1.2mm]
\DrRL p_3(t) =  -1+  (\alpha_1 + u(t) + \mu) p_3(t) 
- u(t) p_4(t) + \beta S(t)( p_1(t) -  p_2(t)),\\[1.2mm]
\DrRL p_4(t) = \left(\alpha_2 + \mu+\delta\right) p_4(t) 
+ \delta (k-1) p_2(t) - \delta k p_3(t),
\end{cases} 
\end{equation}
which is a fractional system of right Riemann--Liouville fractional derivatives, 
$\DrRL$. In turn, the optimality condition of PMP establishes 
that the optimal control is given by
\begin{equation}
\label{eq:ext:cont}
u(t)=\min\left\{\max\left\{0,\dfrac{\left(p_3(t)-p_4(t)\right)
I(t)}{2 B \rho}\right\},u_{\max}\right\}.
\end{equation}
In addition, the following transversality conditions hold:
\begin{equation} 
\label{eq:transversality}
_{\scriptscriptstyle t}D^{\alpha-1}_{\scriptscriptstyle t_f}
p_i\bigm|_{\scriptscriptstyle t_f}=0\Leftrightarrow 
_{\scriptscriptstyle t}\!I^{1-\alpha}_{\scriptscriptstyle t_f}
p_i\bigm|_{\scriptscriptstyle t_f}=p_i(t_f)=0, 
\quad i=1,\ldots,4,
\end{equation}
where $_{\scriptscriptstyle t}\!I^{1-\alpha}_{\scriptscriptstyle t_f}$
is the right Riemann--Liouville fractional integral of order $1-\alpha$.

\subsection{Numerical results and cost-effectiveness of the fractional TB optimal control problem}
\label{sec:nresults}

The PMP is used to numerically solve the optimal control problem 
\eqref{cost-functional}--\eqref{Omega:set}, as discussed in 
Section~\ref{sec:optctrl}, in the classical ($\alpha = 1$)
and fractional ($0 < \alpha < 1$) cases, using the
predict-evaluate-correct-evaluate (PECE) method of Adams--Basforth--Moulton
\cite{diethelm2005algorithms}, implemented in MATLAB. 
First, we solve system \eqref{eq:modTB_control}
by the PECE procedure with initial values for the state variables 
given by  Table~\ref{tab:solinit} and a guess for the control 
over the time interval $[0,t_f]$, and obtain the values 
of the state variables $S$, $L$, $I$ and $T$.
A change of variables is applied to the adjoint system \eqref{eq:co_states_fr2}
and to the transversality conditions \eqref{eq:transversality},
obtaining a left Riemann--Liouville fractional initial value problem. 
Such system is solved with the PECE procedure and
the values of the co-state variables $p_i$, $i=1,\ldots,4$, are obtained.
The control is then updated by a convex combination of the previous control
and the value computed according to \eqref{eq:ext:cont}. This procedure 
is repeated iteratively until the values of the controls at the previous 
iteration are very close to the ones at the current iteration.
The solutions of the classical problem ($\alpha=1$) where successfully 
confirmed by an algorithm that uses a classical forward-backward scheme, 
also implemented in MATLAB.

In our numerical experiments, we consider that the maximum number 
of infectious individuals for the problem without control, $\rho$, 
is $452.758$, $u_{\max}=1$, $B=0.15$, while the other parameters 
are fixed according with Table~\ref{tab:param}. Note that, with 
the above values for the parameters, $R_0 = 7.1343 > 1$ (endemic situation).
Our initial conditions, given by Table~\ref{tab:solinit}, guarantee the existence 
of a non-trivial endemic equilibrium for system \eqref{eq:modTB_control} 
without control. Because the WHO (World Health Organization) goals for most 
diseases are usually fixed for five years periods, we considered $t_f=5$.
\begin{table}[!htb]
\centering
\caption{Initial conditions
for the fractional optimal control problem \eqref{cost-functional}--\eqref{Omega:set}
with parameters given by Table~\ref{tab:param}, corresponding to the endemic equilibrium
of TB model \eqref{TB_model}.}\label{tab:solinit}
\begin{tabular}{c@{\hspace*{2.2cm}}c@{\hspace*{2.2cm}}c@{\hspace*{2.2cm}}c}\toprule
$S(0)$ & $L(0)$ & $I(0)$ & $T(0)$ \\[1mm] \midrule
7779.28 & 43511.9 & 175.267 &  78.4299\\ \bottomrule
\end{tabular}
\end{table}

\begin{figure}[!htb]
\centering
\begin{subfigure}[b]{0.46\textwidth}\centering
\includegraphics[scale=0.46]{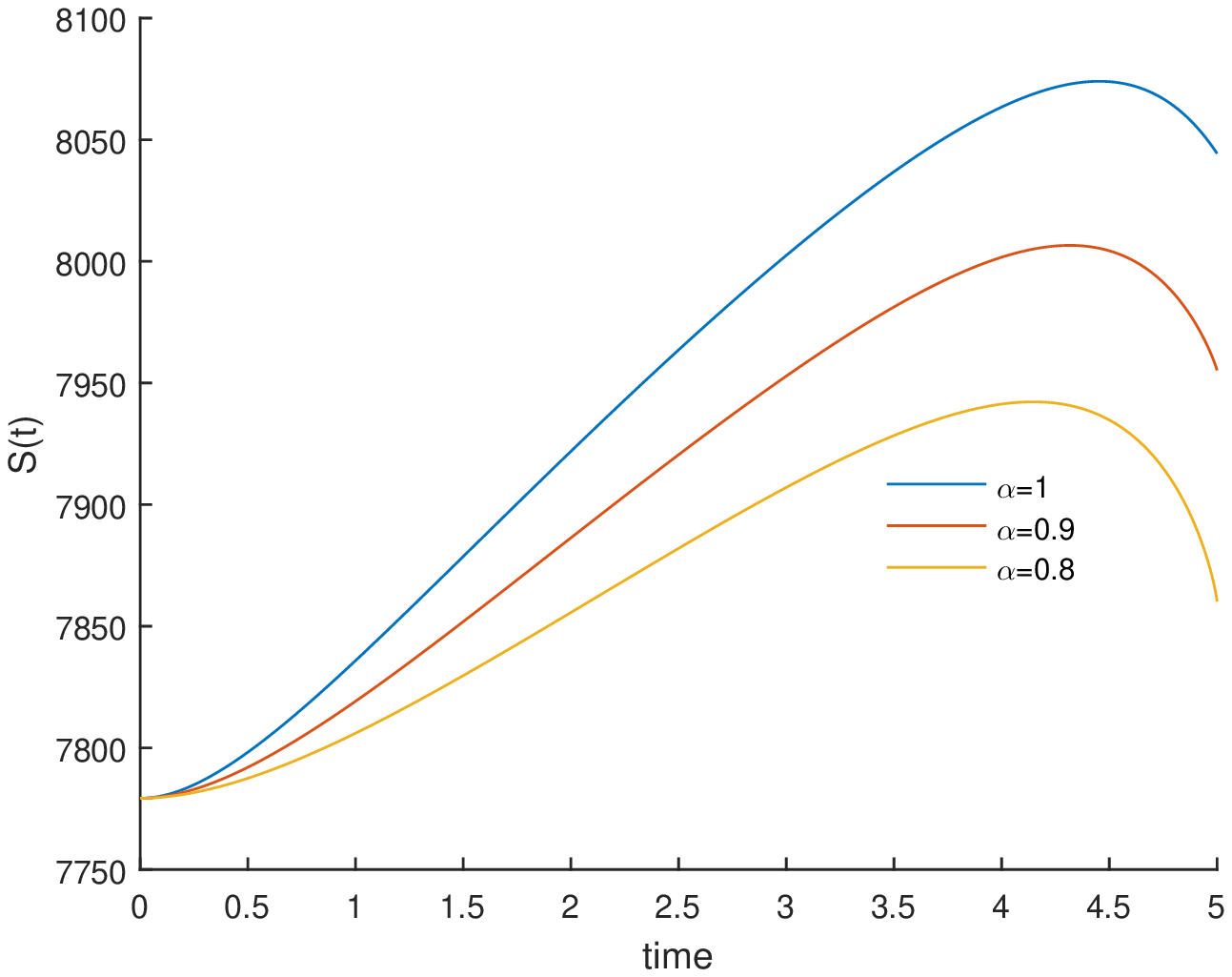}
\caption{Susceptible individuals.}
\end{subfigure}\hspace*{1cm}
\begin{subfigure}[b]{0.46\textwidth}
\centering
\includegraphics[scale=0.46]{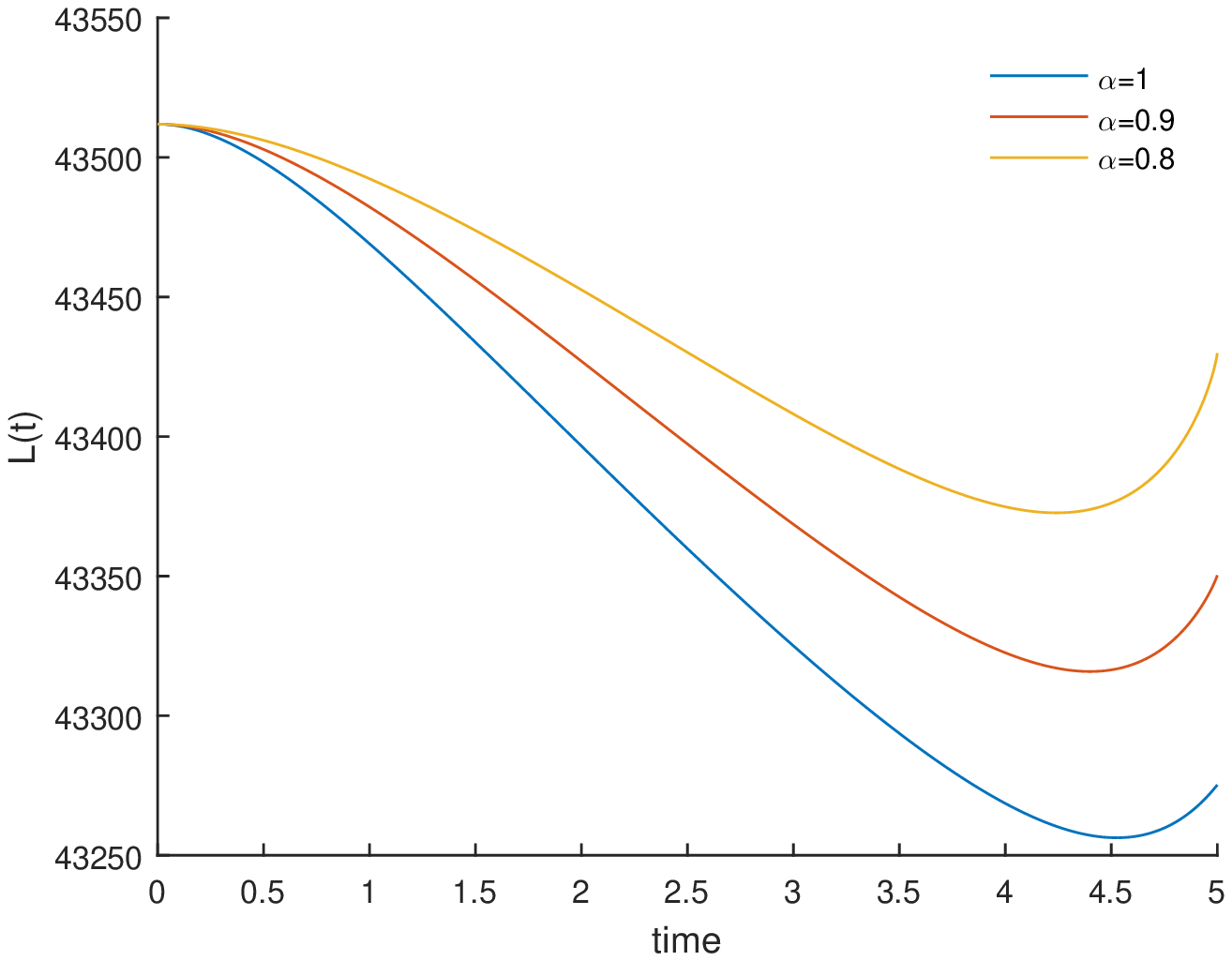}
\caption{Latent individuals.}
\end{subfigure}\\
\begin{subfigure}[b]{0.46\textwidth}
\centering
\includegraphics[scale=0.46]{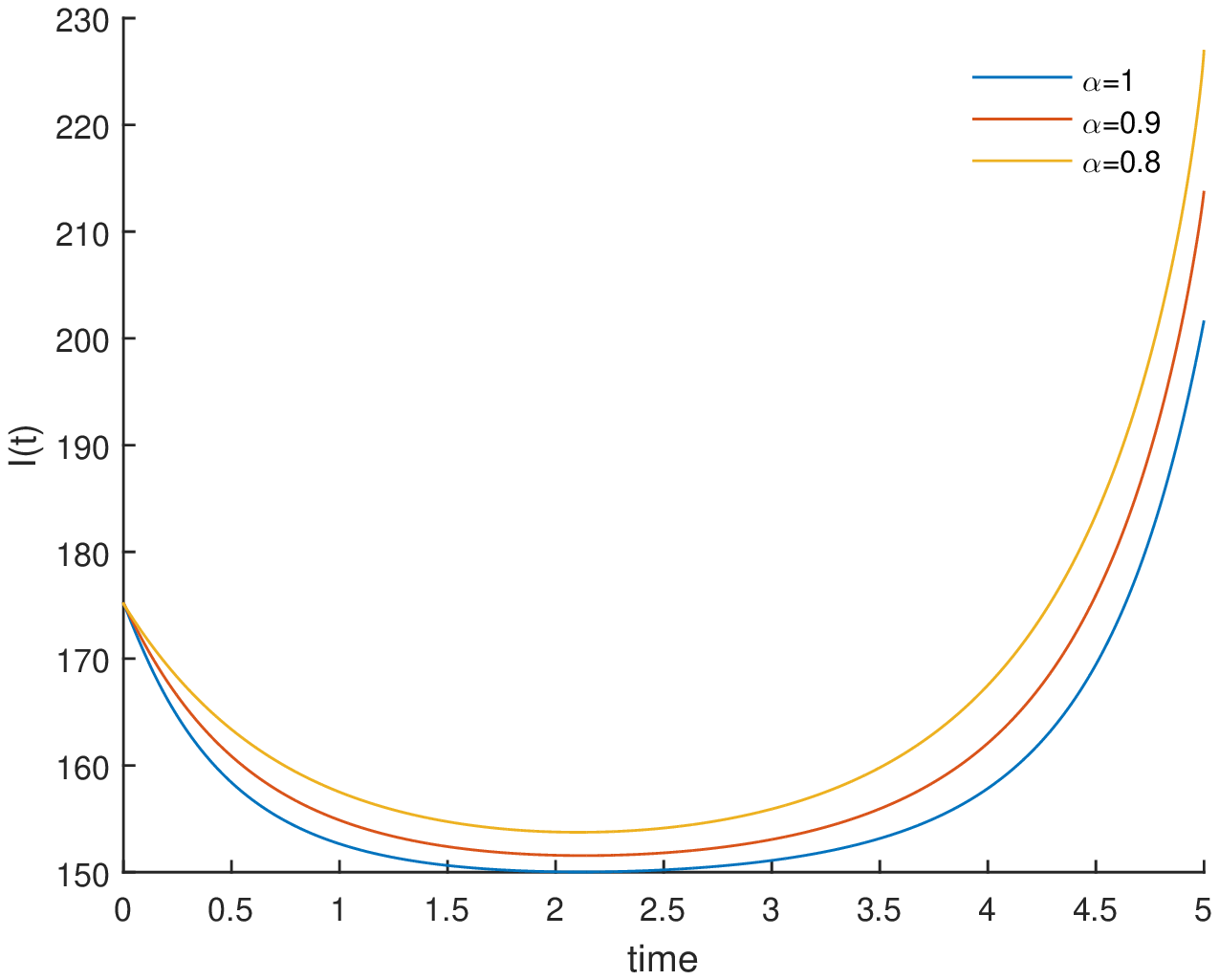}
\caption{Infectious individuals.}\label{fig:I_var:alphas}
\end{subfigure}\hspace*{1cm}
\begin{subfigure}[b]{0.46\textwidth}
\centering
\includegraphics[scale=0.46]{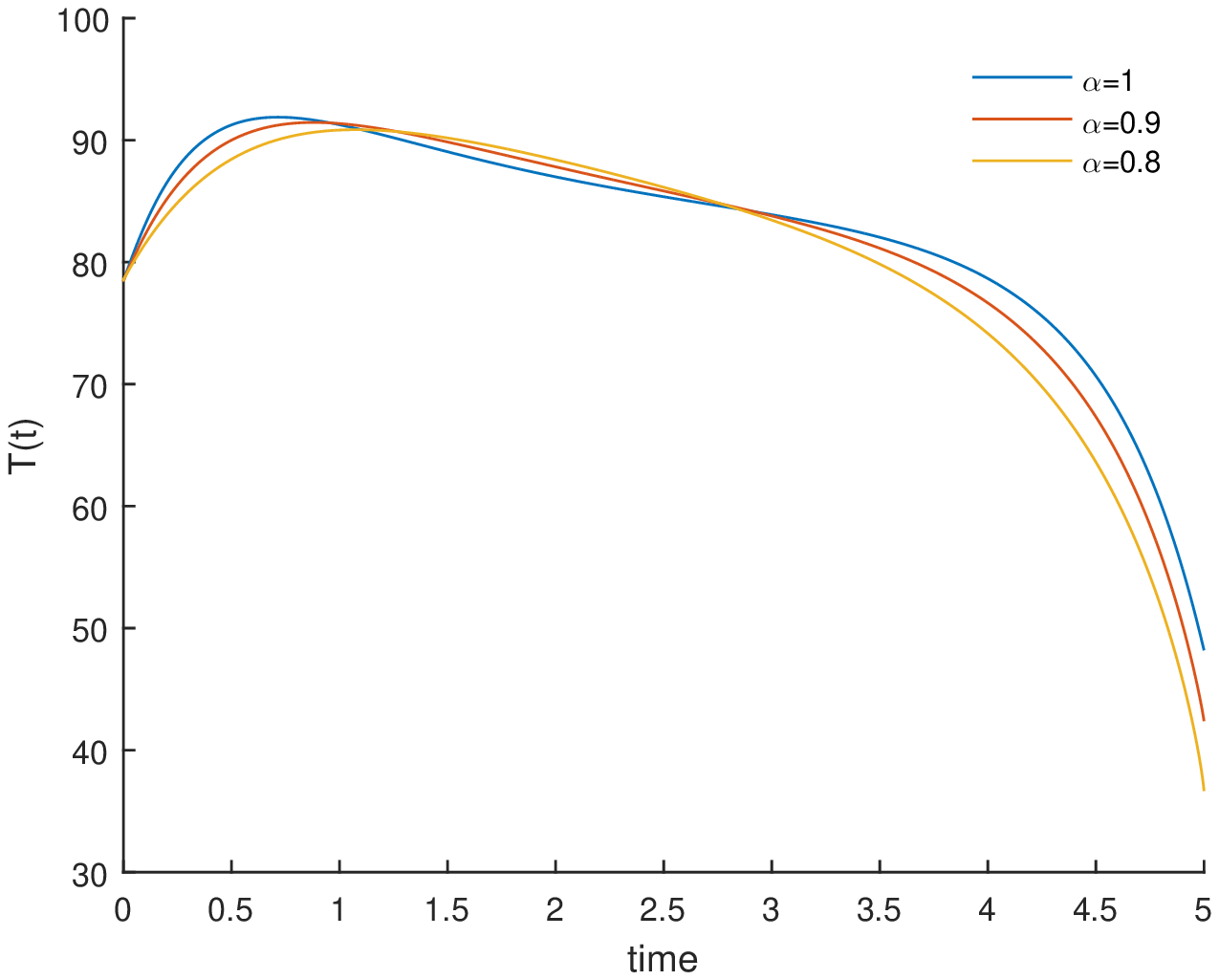}
\caption{Treated individuals.}
\end{subfigure}
\caption{State variables of the fractional optimal control problem \eqref{cost-functional}--\eqref{Omega:set},
with values from Table \ref{tab:param}, weight $B=0.15$, and the fractional order derivatives  $\alpha=1.0,$ $0.9$ and $0.8$.}
\label{fig:states_var:alphas}
\end{figure}

\begin{figure}[!htb]
\centering
\begin{subfigure}[b]{0.46\textwidth}\centering
\includegraphics[scale=0.46]{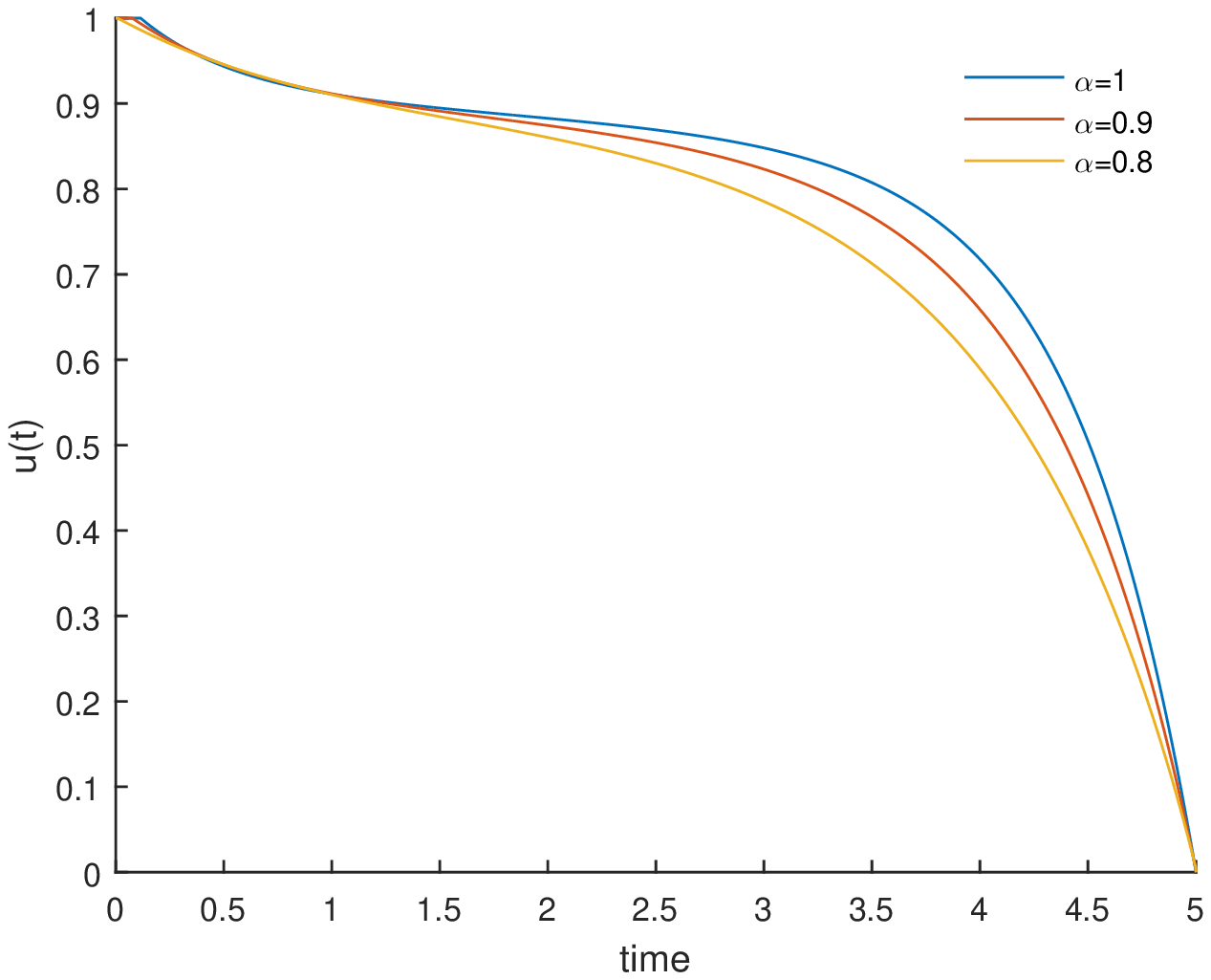}
\caption{Optimal control, u.}\label{fig:u_var:alphas}
\end{subfigure}\hspace*{1cm}
\begin{subfigure}[b]{0.46\textwidth}
\centering
\includegraphics[scale=0.46]{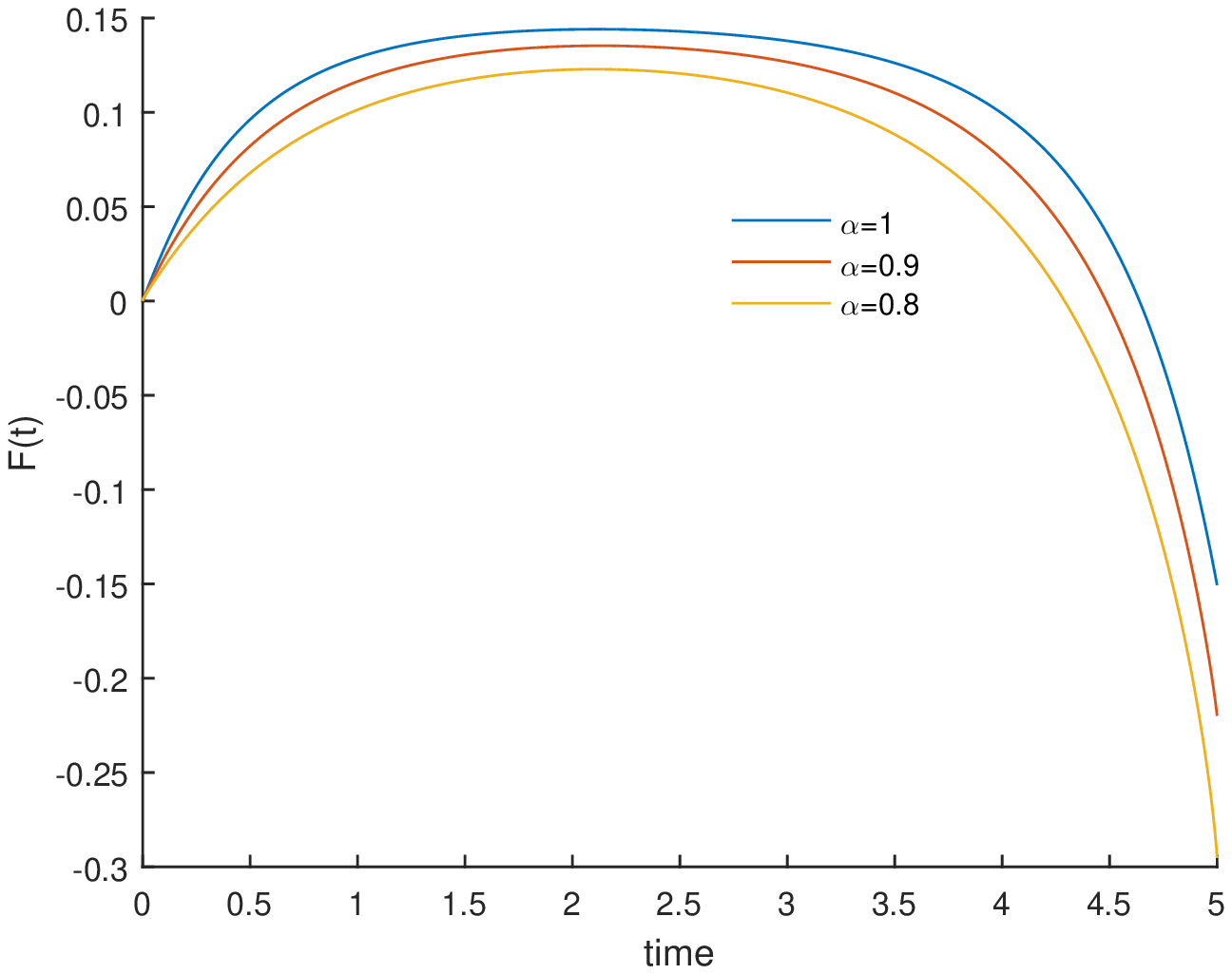}
\caption{Efficacy function $F(t)$.}\label{fig:efficacy:alphas}
\end{subfigure}
\caption{Control function $u(t)$ and efficacy function $F(t)$ associated to the fractional optimal control problem \eqref{cost-functional}--\eqref{Omega:set},
with values from Table \ref{tab:param}, weight $B=0.15$, and the fractional order derivatives  $\alpha=1.0,$ $0.9$ and $0.8$.}
\label{fig:u_F:alphas}
\end{figure}

Without loss of generality, we considered the fractional order derivatives 
$\alpha=1.0,$ $0.9$ and $0.8$. In Figures~\ref{fig:states_var:alphas} 
and \ref{fig:u_var:alphas}, we find the solutions of the fractional optimal 
control problem for those values of $\alpha$. We can see that a change 
in the value of $\alpha$, corresponds significant variations of the state 
and control variables. The existence of an endemic situation ($R_0> 1$), 
a very large initial number of latent individuals (about 84\% of population),  
and a control that vanishes at the end of the time interval motivates that, 
in the end of the time interval, the number of infected individuals  
exceeds its initial value. Beyond those values of $\alpha$, others values 
less than one were also tested, but the results do not change qualitatively.  
We note that, in all the experiments, decreasing $\alpha$ meant increasing 
the number of infected, as Figure~\ref{fig:I_var:alphas} evidences.

In Figure~\ref{fig:efficacy:alphas}, the efficacy function 
\cite{rodrigues2014cost} is displayed. It is defined by
\begin{equation*}
F(t)=\frac{I(0)-I^*(t)}{I(0)}=1-\frac{I^*(t)}{I(0)},
\end{equation*}
where $I^*(t)$ is the optimal solution associated with the fractional 
optimal control and $I(0)$ is the corresponding initial condition. 
This function measures the proportional variation  in the number of
infectious individuals after the application of the control
$u^*$, by comparing the number of infected individuals
at time $t$ with the initial value $I(0)$. We observe that $F(t)$ exhibits
the inverse tendency of $I(t)$. Once $I(t)$ ends with values bigger 
than its initial value,  $F(t)$ ends with negative values.

To assess the cost and the effectiveness of the proposed fractional
control measure during the intervention period,
some summary measures are presented.
The total cases averted by the intervention, during
the time period $t_f$, is defined in \cite{rodrigues2014cost} by
\begin{equation*}
A=t_f I(0)-\int_0^{t_f}I^*(t)~dt,
\end{equation*}
where $I^*(t)$ is the optimal solution associated with the fractional 
optimal control $u^*$ and $I(0)$ is the corresponding initial condition. 
Note that this initial condition is obtained as the equilibrium proportion
$\overline{I}$ of system \eqref{eq:modTB_control} without treatment
intervention, which is independent on time, so that
$t_f I(0)=\displaystyle \int_0^{t_f}\overline{I}~dt$ represents the total 
infectious cases over the given period of $t_f$ years.
We define effectiveness as the proportion of cases averted
on the total cases possible under no intervention \cite{rodrigues2014cost}:
\begin{equation*}
\overline{F}=\frac{A}{t_f I(0)}
=1-\frac{\displaystyle \int_0^{t_f}I^*(t)~dt}{t_f I(0)}.
\end{equation*}
The total cost associated with the intervention
is defined in \cite{rodrigues2014cost} by
\begin{equation*}
TC=\int_0^{t_f} C \, u^*(t)I^*(t)~dt,
\end{equation*}
where $C$ corresponds to the unit cost, per person, of detection
and treatment of infectious individuals.
Following \cite{okosun2013optimal,rodrigues2014cost},
the average cost-effectiveness ratio is given by
\begin{equation*}
ACER=\frac{TC}{A}.
\end{equation*}
In Table~\ref{tab:efficacy}, the cost-effectiveness measures, 
for the fractional optimal control problem
we have analysed, are summarized. These results
show the effectiveness of the control to reduce TB 
infectious individuals and the superiority of the classical model ($\alpha=1$).
\begin{table}[!htb]
\centering
\caption{Summary of cost-effectiveness measures for classical and fractional 
($0<\alpha <1$) TB disease optimal control problems. Parameters according
to Tables~\ref{tab:param} and~\ref{tab:solinit} with $C=1$.}\label{tab:efficacy}
\begin{tabular}{c@{\hspace*{1cm}}c@{\hspace*{1cm}}c@{\hspace*{1cm}}c@{\hspace*{1cm}}c} 
\toprule
$\alpha$ & $A$  & $TC$ & $ACER$ & $\overline{F}$  \\[1mm] \midrule
1.0   &  88.6409  &  617.147  &   6.96233  &  0.101150\\
0.9   & 73.4419   & 608.801   &  8.28956   & 0.083806\\
0.80  &  54.8307  &  597.161  &   10.891   & 0.062568\\
 \bottomrule
\end{tabular}
\end{table}

\begin{figure}[!htb]
\centering
\begin{subfigure}[b]{0.46\textwidth}\centering
\includegraphics[scale=0.46]{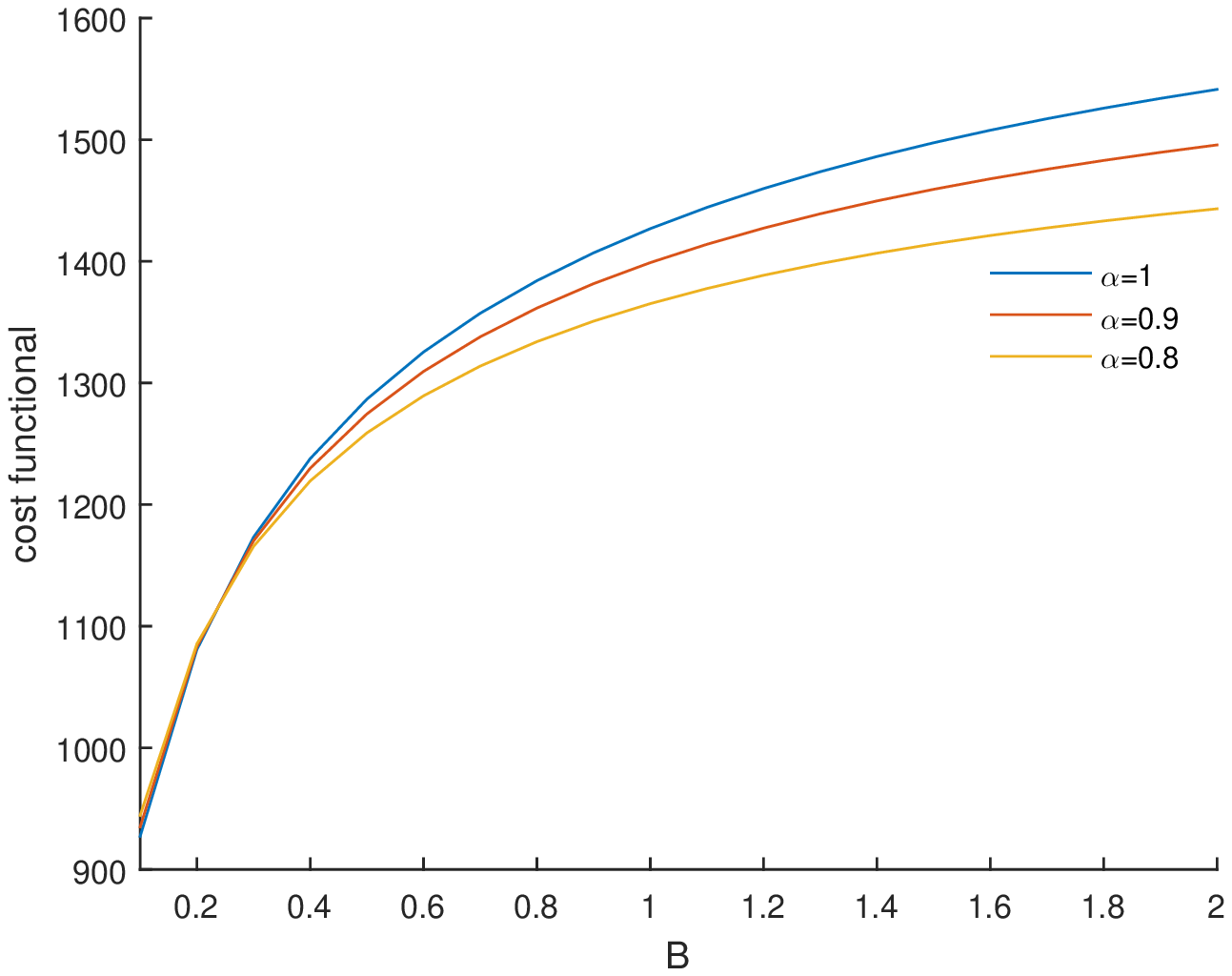}
\caption{Evolution of the cost functional $\mathcal{J}$.}\label{fig:evol_J}
\end{subfigure}\hspace*{1cm}
\begin{subfigure}[b]{0.46\textwidth}
\centering
\includegraphics[scale=0.46]{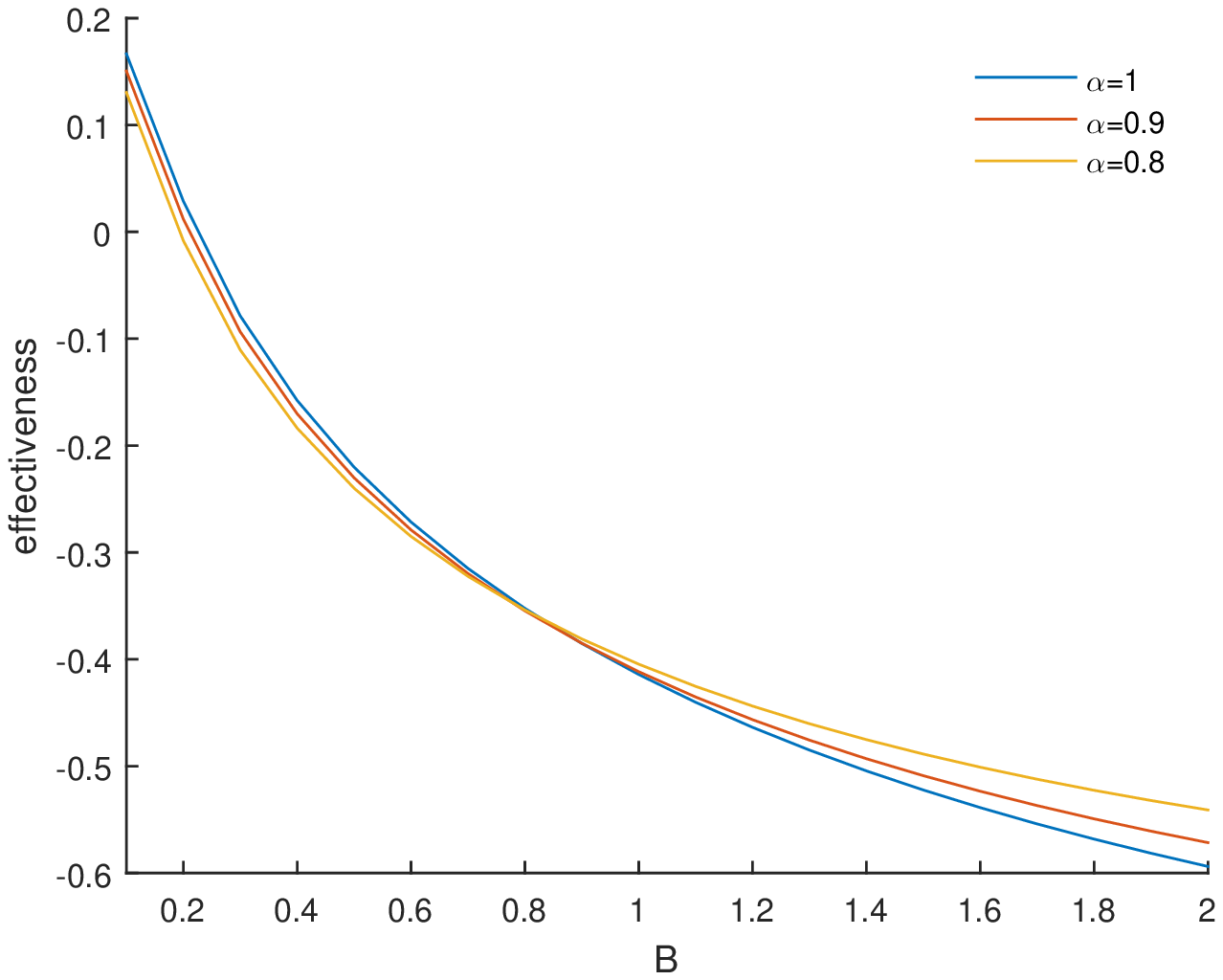}
\caption{Evolution of the effectiveness $\overline F$.}\label{fig:evol_F}
\end{subfigure}
\caption{Impact of the variation of the weight $B$  on the cost functional value, 
$\mathcal{J}$, (left) and on the effectiveness measure $\overline{F}$ (right) 
for fractional order derivatives $\alpha=1.0$, $0.9$ and $0.8$.}
\label{fig:evol_J_F}
\end{figure}

The impact of the variation of the weight $B$ on 
the cost functional is displayed in Figure~\ref{fig:evol_J}, 
and the effectiveness measure $\overline{F}$ is displayed 
in Figure~\ref{fig:evol_F}. When the cost of treatment increases, 
we observe that the fractional model can be more effective in reducing 
the number of infective individuals at a not so high cost 
(Figure~\ref{fig:evol_J}).

\section{Conclusions}
\label{sec:conclusions}

A sensitivity analysis was carried out for the basic
reproduction number of the TB model \cite{TB:frac:2018}, 
proving that the most important parameter to have into account is
the natural death rate $\mu$. The application of optimal control 
to the fractional TB model shows that treatment is effective
on the reduction of infected individuals. A change in the value 
of the fractional derivative order, $\alpha$,
corresponds to substantial variations on the solutions 
of the fractional optimal control problem.
The fractional model ($0<\alpha<1$) is recommended
when treatment is expensive, and only in that case, 
because in that situation it is more effective 
and not so costly as the classical model.


\section*{Acknowledgements}

Rosa was supported by the Portuguese Foundation for Science and Technology (FCT)
through IT (project UID/EEA/50008/2013);
Torres by FCT through CIDMA (project UID/MAT/04106/2019)
and TOCCATA (project PTDC/EEI-AUT/2933/2014
funded by FEDER and COMPETE 2020).



\end{document}